\documentclass[journal]{IEEEtran}

\usepackage[pdftex]{graphicx}
\usepackage[cmex10]{amsmath}
\usepackage{algorithm}
\usepackage{algorithmic}
\usepackage{array}
\usepackage{fixltx2e}
\usepackage{amsfonts}
\usepackage{pgfplots}
\usepackage{hyperref}  
\usepackage{url}  

\author{Michael~Azzam, S.~Easter~Selvan, Augustin~Lef\`evre, and P.-A.~Absil}
\title{Mixed Integer Programming to Globally Minimize the Economic Load Dispatch Problem With Valve-Point Effect}

\newcommand{\R}{\mathbb{R}}
\renewcommand{\vec}[1]{\mathbf{#1}}
\newtheorem{theorem}{Theorem}
\newtheorem{proposition}[theorem]{Proposition} 

\begin{document}
\maketitle

\begin{abstract}

Optimal distribution of power among generating units to meet a specific demand subject to system constraints is an ongoing research topic in the power system community. The problem, even in a static setting, turns out to be hard to solve with conventional optimization methods owing to the consideration of \emph{valve-point effects} which make the cost function nonsmooth and nonconvex. This difficulty gave rise to the proliferation of population-based global heuristics in order to address the multi-extremal and nonsmooth problem. In this paper, we address the economic load dispatch problem (ELDP) with valve-point effect in its classic formulation where the cost function for each generator is expressed as the sum of a quadratic term and a rectified sine term. We propose two methods that resort to piecewise-quadratic surrogate cost functions, yielding surrogate problems that can be handled by mixed-integer quadratic programming (MIQP) solvers. The first method shows that the global solution of the ELDP can often be found by using a fixed and very limited number of quadratic pieces in the surrogate cost function. The second method adaptively builds piecewise-quadratic surrogate under-estimations of the ELDP cost function, yielding a sequence of surrogate MIQP problems. It is shown that any limit point of the sequence of MIQP solutions is a global solution of the ELDP. Moreover, numerical experiments indicate that the proposed methods outclass the state-of-the-art algorithms in terms of minimization value and computation time on practical instances.
\end{abstract}


\section{Introduction}
\label{sec:introduction}

\IEEEPARstart{E}{conomic load dispatch problem} (ELDP) attempts to minimize the cost associated with the power generation by optimally scheduling the load across generating units to meet a certain demand subject to system constraints~\cite{Sinha03}. It is not uncommon to notice that the cost function for a generator is approximated by a quadratic function for the sake of simplicity. Nevertheless, when the cost function also takes into account highly nonlinear input-output characteristics due to valve-point loadings, even the static ELDP problem that ignores the ramp-rate constraints turns out to be difficult to solve. The challenges faced by the solvers stem from (i) the nonsmooth cost function and (ii) the multi-extremal nature of the problem.

\par A popular strategy to address the ELDP is to rely on a (population-based) stochastic search algorithm. Indeed, a myriad of such algorithms have been proposed during recent years, including genetic algorithms~\cite{Walters93}, evolutionary programming~\cite{Sinha03}, particle swarm optimization~\cite{Park05}, ant swarm optimization~\cite{CaiAug07}, differential evolution~\cite{Coelho06}, firefly algorithm~\cite{Yang11}, bacterial foraging algorithm~\cite{Panigrahi08}, and biogeography-based optimization~\cite{BhattacharyaNov10}. 

\par Since heuristics enable the global exploration, and a local method aids to converge to a local optimum, a focussed effort has been made by many to integrate a global heuristics with a local optimizer, resulting in hybrid algorithms. Interested readers may refer to~\cite{Subathra14} for an exhaustive list of such algorithms in the ELDP context. A few well-known local optimizers integrated with a global scheme to tackle the ELDP are the Nelder-Mead method~\cite{Panigrahi08}, generalized pattern search~\cite{Alsumait10}, sequential quadratic programming~\cite{Victoire04}, and a Riemannian subgradient steepest descent~\cite{Borckmans14}. Note that the equality constraint must be handled by way of a slack variable or a barrier approach, and the inequality constraints with the help of a penalty approach in global heuristics. Furthermore, even though the global heuristics favor finding the global minimum, these hybrid algorithms are guaranteed at best to find a local minimum.

In this paper, we introduce new techniques to efficiently build surrogates of the ELDP cost function in order to take advantage of powerful modern mixed-integer programming (MIP) solvers. In particular, the adaptive technique introduced in Section~\ref{sec:global_convergence} builds a sequence of surrogate piecewise-quadratic cost functions that aims at keeping low the number of pieces for the sake of efficiency while nevertheless offering the guarantee that the sequence of surrogate solutions converge to the global solution of the ELDP. It is interesting to note that the MIP method vested with theoretical guarantees as applied to a static ELDP with the valve-point effect may be extended to a dynamic setting, provided the number of generating units is not too large. Finally, we demonstrate that the minimization results from a 3-, two instances of 13-, and a 40-generator settings are lower than the results reported thus far in the literature with the same datasets.          

While the rectified-sine model~\eqref{eq:fi} of the valve-point effect was proposed more than 20 years ago~\cite{Walters93}, it is only in the past few months that MIP techniques appeared in the literature to handle this problem formulation~\cite{Yang13,Pedroso14,Wang14}. The methods introduced in this paper contribute beyond this recent literature in two ways. (i) The approach proposed in Section~\ref{sec:simple_solver} shows that it is often possible to obtain the exact global solution with a fixed and very limited number of linear pieces in the surrogate cost function. (ii) The adaptive approach proposed in Section~\ref{sec:global_convergence}, while sharing several aspects with the recent report~\cite{Pedroso14}, introduces fewer breakpoints in the surrogate cost function, allowing for a reduced complexity. Moreover, it takes into account that it is only the rectified sine term that is piecewise concave in~\eqref{eq:fi}; this leads a piecewise-quadratic under-approximation of~\eqref{eq:fi} that is handled by mixed-integer quadratic programming (MIQP).

The remainder of the article is organized as follows. In Section~\ref{sec:problem}, the ELDP with the valve-point effect is briefly presented. The principle behind our approach is first expounded in Section~\ref{sec:simple_solver} using a simple linear approximation of the term that accounts for the nonconvexity and nonsmoothness. Section~\ref{sec:finer_solver} reviews MIP formulations for general piecewise-linear objective functions and shows how this technique can be integrated in the ELDP.
The adaptive algorithm is introduced and its convergence analysis carried out in Section~\ref{sec:global_convergence}. Numerical results are reported in Section~\ref{sec:results}, and conclusions are drawn in Section~\ref{sec:conclusion}.

\section{Problem statement}
\label{sec:problem}

In this section, we recall the formulation of the widely investigated ELDP with valve-point effect, as described, e.g., in~\cite{Subathra14}.

\begin{figure}
\centering
\begin{tikzpicture}[scale=1]
\begin{axis}[ 
	axis x line=bottom,
	axis y line=left,
	xmin=30,xmax=220,
	ymin=800,ymax=2700,
	xlabel=$p$, ylabel=$f(p)$
	]
	\addplot[id=1a,domain=50:200] gnuplot[samples=250] {0.00533*x*x+11.669*x+213.1};
	\addplot[id=1b,domain=50:200,dashed] gnuplot[samples=250] {0.00533*x*x+11.669*x+213.1+130*abs(sin(0.0635*(x-50)))};
\end{axis}
\end{tikzpicture}
\caption{Examples of cost functions for a generator, without (solid) or
with (dashed) valve-point effect}
\label{fig:example_function}
\end{figure}
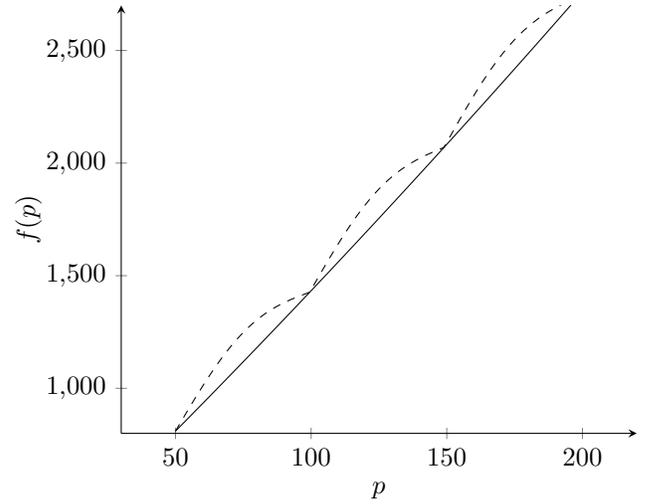

In the ELDP, the main component of the cost that needs to be taken into account, is the cost of the input, that is to say the fuel. The objective function used
in the problem is thus defined by how we represent the input-output
relationship of each generator. The total cost is then naturally the sum
of each contribution. The objective function is thus written
\begin{subequations}  \label{eq:f}
\begin{equation}
f(\vec{p}) = \min_{\vec{p}\in \R^n} \sum_{i=1}^n f_i(p_i),
\end{equation}
where $f$ is the total cost function in \$/h, equal to the sum of the $n$
$f_i$ univariate functions that give the individual contribution in the
total cost of the $i$th generator, depending on the $p_i$ amount of power, in MW, assigned to this unit.

A classical, simple and straightforward approach to construct the cost
functions is to use a quadratic function for each generator, i.e.,
$
 f_i(p_i) = a_i p_i^2 + b_i p_i + c_i,
$
where $a_i$, $b_i$ and $c_i$ scalar coefficients.

However, in reality, performance curves do not behave so smoothly. In the
case of generating units with multi-valve steam turbines, ripples will
typically be seen in the curve. Large steam turbine generators usually have a number of steam admission valves that are opened in sequence to meet an increasing demand from a unit. And as each steam admission
valve starts to open, a sharp increase in losses due to wire drawing effects occur~\cite{Subathra14}, \cite{IEEE_report}. This is the so-called valve-point effect.
To try to capture this effect in the model, a rectified sine term is usually added
to the fuel cost functions, so that they become
\begin{equation}  \label{eq:fi}
 f_i(p_i) = a_i p_i^2 + b_i p_i + c_i + d_i |\sin(e_i(p_i-p^{min}_i))|
\end{equation}
\end{subequations}
where $d_i$ and $e_i$ are additional positive coefficients needed to
take the valve-point effect into account~\cite{Walters93}.
We can see how the new term affects the cost function in the example of \figurename~\ref{fig:example_function}.

Unfortunately, this addition brings two detrimental properties in the
problem: non-convexity and non-differentiability. These are two
major hindrances that prevent the direct application of usual
optimization algorithms.

Naturally, the problem has also some constraints that must be satisfied
and restrict the search space. First, the producer must meet the demand,
even though some power will be lost in the network. This is the \emph{power balance constraint}, formulated through an equality constraint
\begin{equation}  \label{eq:power-balance}
\sum_{i=1}^n p_i = p_D + p_L(\vec p)
\end{equation}
with scalar $p_D$ and function $p_L$ being respectively the demand on the system and the
transmission loss in the network, both in MW. 
The transmission loss is computed using the so-called \emph{B-coefficients} as
\begin{equation}
p_L(\vec p) = \vec{p}^t B \vec{p} + \vec{p}^t \vec b^0 + b^{00}
\end{equation}
where $B$ is a symmetric positive-semidefinite matrix, $b^0$ a vector of
size $n$ and $b^{00}$ a scalar.

The other type of constraints are \emph{generator capacity constraints}, which take the form of box constraints imposing that each unit
has its own range of possible power generation, from $p^{min}_i$ to
$p^{\mathrm{max}}_i$. This is easily transcribed as inequality constraints
\begin{equation}  \label{eq:box}
p^{\mathrm{min}}_i \leq p_i \leq p^{\mathrm{max}}_i,
\end{equation}
where $p^{\mathrm{min}}_i$ and $p^{\mathrm{max}}_i$ are obviously the minimum and maximum
power output of the $i$th generator, in MW.

The set of points that satisfy constraints~\eqref{eq:power-balance} and~\eqref{eq:box} is termed the \emph{feasible set} of the ELDP.

As is often done, we will
ignore the transmission loss in the network, i.e., we set $p_L = 0$ in~\eqref{eq:power-balance}. 
The static ELDP without losses is then classically written as
\begin{equation}
\begin{aligned}
\min_{\vec p \in \R^n} ~&\sum_{i=1}^n a_ip_i^2+b_ip_i+c_i + d_i|\sin(e_i(p_i-p^{\mathrm{min}}_i))| \\
\textrm{subject to} ~~ &\sum_{i=1}^n p_i = p_D \\
    &p^{\mathrm{min}}_i \leq p_i \leq p^{\mathrm{max}}_i \quad \forall i \in \{1,...,n\}
\end{aligned}
\label{eq:eldp}
\end{equation}
However, as we point out in the concluding section, considering the loss does not add much complexity. 

\section{A first simple MIQP approach}
\label{sec:simple_solver}

We now introduce our first, crude but effective way of building a piecewise-quadratic surrogate of the ELDP cost function~\eqref{eq:f}, and we show how to solve the resulting surrogate problem with an MIQP solver.

As discussed in the previous section, the term 
\begin{equation}
|\sin(e_i(p_i-p^{\mathrm{min}}_i))|
\label{sin}
\end{equation}
of the objective function makes the optimization problem challenging because it breaks its smoothness and convexity.

In this paper, we overcome the difficulty by approximating~\eqref{sin} with functions that are more manageable, namely piecewise-linear functions. Even though they do not restore smoothness nor convexity of the problem, they are conveniently handled by mixed-integer programming.

A first simple piecewise-linear approximation consists of replacing $|\sin x|$ by $|x|$ over $[-\pi/2,\pi/2]$ and completing the approximation over the whole domain by periodic extension, observing that $|\sin x|$ is periodic of period $\pi$. The underlying motivation is to keep low the number of linear pieces while capturing accurately the behavior of~\eqref{sin} around its kink points, as the optimum tends to be located at those kink points.

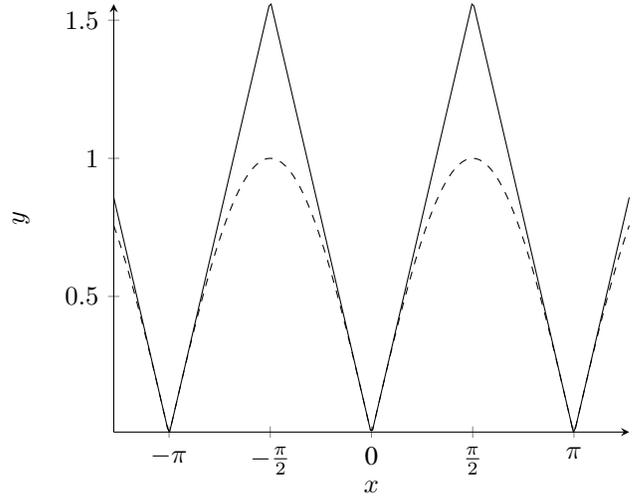
\begin{figure}
\centering
\begin{tikzpicture}
\begin{axis} [
	axis x line=bottom,
	axis y line=left,
	xlabel=$x$, ylabel=$y$,
	xtick={-3.14159, -1.5708, 0, 1.5708, 3.14159},
	xticklabels={$-\pi$, $-\frac{\pi}{2}$, 0, $\frac{\pi}{2}$, $\pi$}
	]
  \addplot[id=2a,domain=-4:4,dashed] gnuplot[samples=250] {abs(sin(x))};
  \addplot[id=2b,domain=-4:4,solid] gnuplot[samples=250]{abs(asin(sin(x)))};
\end{axis}
\end{tikzpicture}

\caption{Plot of $|\sin(x)|$ on the interval [-4,4] and the suggested surrogate : $|\arcsin(\sin(x))|$. Sine is in dashed line and the surrogate is solid.}
\label{approx}

\end{figure}

The resulting function can be compactly written as $|\arcsin(\sin x)|$, which has the sawtooth shape shown on \figurename~\ref{approx}. This function can be interpreted as the distance between $x$ and the closest
multiple of $\pi$. This can be written as a small mixed integer optimization problem such as
\[ |\arcsin(\sin x)| = \min_{k\in\mathbb{Z}} |x-k\pi| \]
which can be reformulated to become linear, as follows :
\begin{align*}
\min_{(k,t) \in \mathbb{Z}\times\R}~ & t\\
\text{s.t. } ~ & x-k\pi \leq t \\
               & x-k\pi \geq -t .
\end{align*}

The problem is thus very simple, even though we introduce an integer variable. By introducing this in the original problem, we get
\begin{equation}
\begin{aligned}
\min_{p,k,t} &\sum_i a_ip_i^2+b_ip_i+c_i + d_i t_i \\
\mathrm{s.t. } ~ &\sum_i p_i = D \\
    &p^{\mathrm{min}}_i \leq p_i \leq p^{\mathrm{max}}_i \\
    & -t_i \leq e_i(p_i-p^{\mathrm{min}}_i)-\pi k_i \leq t_i \\
    & p_i \in \mathbb{R} \\
    & t_i \in \mathbb{R} \\
    & k_i \in \mathbb{N} \\
\end{aligned}
\label{eq:simple}
\end{equation}
which is an MIQP problem. This
class of problems can be solved exactly by solvers, for instance, CPLEX, Gurobi or Mosek. In Section~\ref{sec:results}, we see that this model
gives results that are competitive with other methods suggested in the literature.

\section{Finer approximation with piecewise linear functions}
\label{sec:finer_solver}

In this section, we show how to handle general piecewise-linear surrogate objective functions, then we consider specifically an over-approximation obtained from tangents to the rectified sine function. 

Let us say we want to minimize on an interval $[v, w]$, a piecewise linear function $g$ described by the slopes $\alpha_i$ and intercepts $\beta_i$ of its $m$ line segments components as well as its breakpoints $X_1=v,X_2,\dots, X_{m+1}=w$. This objective function can be expressed as a mixed integer linear problem such as

\begin{equation}
\begin{aligned}
\min ~ &\sum_{i=1}^m \alpha_i \chi_i + \beta_i \eta_i \\
\text{s.t. }~ &\sum_{i=1}^m \chi_i = x \\
&\sum_{i=1}^m \eta_i = 1 \\
&X_i \eta_i  \leq \chi_i \leq X_{i+1} \eta_i  \quad \forall i \in \{1,\dots, m\} \\
&\eta_i \in \{0,1\}  \quad \forall i \in \{1,\dots, m\}.
\end{aligned}
\end{equation}

This system of constraints ensures that for a given $x$, only one of
the binary variables $\eta_i$ will be equal to one, indicating which 
segment is active, while one of the real variables $\chi_i$ will hold
the value of $x$ and the others will be equal to zero. Thus, only
one term of the objective will be non zero. It will be equal to the value
of the linear function of slope $\alpha_i$ and intercept $\beta_i$, at this abscissa.

We now integrate this technique in our model. Instead of replacing 
$|\sin x|$ by $\min_{k\in\mathbb{Z}} |x-k\pi|$, we use the output of this
as the abscissa for finer approximation through piecewise linear functions. 
For example, with $m_i$ segments for the cost function of unit $i$,

\begin{equation}
\label{eq:sys_algo}
\begin{aligned}
\min_{p,k,t,\chi,\eta} 
&\sum_i \left[ a_i p_i^2+b_ip_i+c_i+ d_i \sum_{j=1}^{m_i} (\alpha_{ij} \chi_{ij} + \beta_{ij} \eta_{ij} ) \right] \\
s.t.~~ &\sum_i p_i = D \\
    &p^{\mathrm{min}}_i \leq p_i \leq p^{\mathrm{max}}_i \\
    & -t_i \leq e_i(p_i-p^{\mathrm{min}}_i)-\pi k_i \leq t_i \\
    & X_{j} \eta_{i,j} \leq \chi_{i,j} \leq X_{j+1} \eta_{i,j} \\
    & t_i = \sum_{j=1}^{m_i} \chi_{i,j} \\
    & \sum_{j=1}^{m_i} \eta_{i,j} = 1 \\
    & p_i \in \mathbb{R} \\
    & t_i,\chi_{i,j} \in \mathbb{R}^+ \\
    & k_i \in \mathbb{N} \\
    & \eta_{i,j} \in \{0,1\} \\
\end{aligned}\hspace{-20pt}
\end{equation}
where $t_i$ are the distances to the nearest root, $\chi_{i,j}$ and $\eta_{i,j}$ are auxiliary variables that help us get the right value in the objective function depending on which segment we are located in and
$\alpha_{ij}$ and $\beta_{ij}$ are constants that parametrize the line
segments.

A course of action we could opt for is following the logic of the 
previous section and build an over-approximation with first-order
Taylor approximation to the sine at different points. 

With three segments, it would mean we want to approximate $\sin x$ by 
$\min(x, T_1(x), T_2(x))$, where $T_1(x)$ and $T_2(x)$ are the 
tangents to $\sin x$ at respectively $\theta_1$ and $\theta_2$.
A bit of algebra  gives us the parameters  \[ \alpha_j = \cos\theta_j \] \[\beta_j= \sin\theta_j-\theta_j\cos\theta_j\] and the endpoints of 
the intervals
\[ X_1 = \frac{\sin\theta_1-\theta_1\cos\theta_1}{1-\cos\theta_1} \]
\[ X_2 = \frac{\theta_2\cos\theta_2-\theta_1\cos\theta_1-\sin\theta_2+\sin\theta_1}{\cos\theta_2-\cos\theta_1} \]

\begin{figure}
\centering
\begin{tikzpicture}
\begin{axis} [
	axis x line=bottom,
	axis y line=left,
	xlabel=$x$, ylabel=$y$,,
	xtick={-3.14159, -1.5708, 0, 1.5708, 3.14159},
	xticklabels={$-\pi$, $-\frac{\pi}{2}$, $0$,, $\pi$},
	extra x ticks={1.0996,1.4765},
 	extra x tick labels={$\theta_1$,$\theta_2$},
    extra x tick style={text height=\heightof{0},},
	]
]
  \coordinate (a) at (axis cs:1.0996,0);
  \coordinate (b) at (axis cs:1.0996,0.7868);
  
  \addplot[id=3a,domain=-4:4,dashed] gnuplot[samples=250] {abs(sin(x))};
  \addplot[id=3b,domain=-4:4] gnuplot[samples=250]
  						 {-0.2730*abs(abs(asin(sin(x)))-0.7176)
  						  -0.1799*abs(abs(asin(sin(x)))-1.2915)
  						  +0.5471*abs(asin(sin(x)))+0.4283};
  \draw[dotted] (axis cs:1.0996,0) -- (axis cs:1.0996,0.8910);
  \draw[dotted] (axis cs:1.4765,0) -- (axis cs:1.4765,0.9956);
\end{axis}
\end{tikzpicture}

\caption{Plot of $|\sin(x)|$ on the interval [-4,4] and a piecewise linear approximation, built from tangents. Sine is in dashed line and the approximation is solid. $\theta_1$ and $\theta_2$ are the tangency
points.}
\label{approx2}
\end{figure}
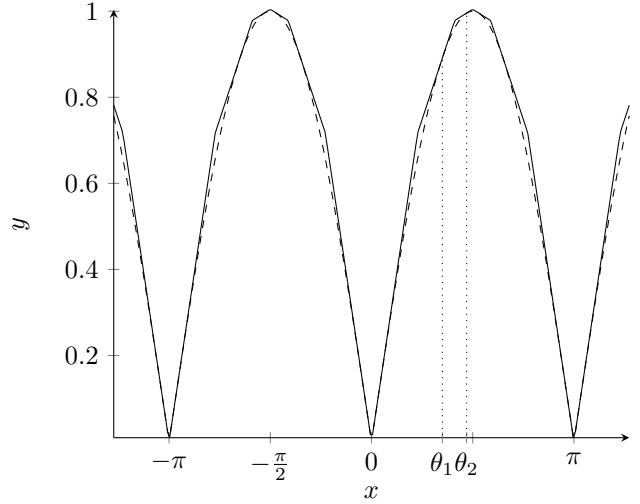

This is illustrated in \figurename~\ref{approx2}. Numerical results are reported in Section~\ref{sec:results}.

\section{Global method based on adaptive under-approximation}
\label{sec:global_convergence}

In this section, we will describe an algorithm that provably converges
toward a global minimum of the ELDP. The way to achieve this is to change our
approach to use under-approximations instead of over-approximations. 

Over-approximations with tangents seem natural because we can then 
infer properties from Taylor's theorem, but under-approximations have the useful property that if the under-approximation and the true function coincide at a global minimizer $x^*$ of the under-approximation, then $x^*$ is also a global minimizer of the true function. 

\begin{proposition}
Let $f:X\mapsto\R$ and $g:X\mapsto\R$ be two functions such that
$g(x) \leq f(x) ~ \forall x\in X$ and $x^* \in X$ a global minimizer
of $g$. If $g(x^*) = f(x^*)$, then $x^*$ is a global minimizer of $f$.
\end{proposition}
\begin{IEEEproof}
Since $x^*$ is a global minimizer of $g$, $g(x^*) \leq g(x)$ for every
$x \in X$. And because $g \leq f$, we have $g(x^*) \leq g(x) \leq f(x)$
for every $x \in X$. Finally, since $f(x^*)=g(x^*)$, we can conclude.
\end{IEEEproof}

So if we were able to find such an under-approximation, then we would have found
the global minimum. What we suggest is to build a sequence of
under-approximations that achieve this goal at the limit.

As before, we will make use of piecewise linear functions to approach the
sine part of the real cost function.
We start off with a simple chord of $|\sin x|$ that links its extreme points. We then solve problem~\eqref{eq:sys_algo}.
We want the next approximation in the sequence to be equal to the true
objective function at the solution $\vec p^0$ we found for the first approximation.
Therefore, we use the values of $t_i$ as breaking points for the new piecewise linear approximation and compute the coefficients of the line
segments so that $g^1(\vec p^0)=f(\vec p^0)$.  We can now simply repeat
this procedure again, find the new solution $\vec p^1$, build a new 
approximation $g^2$, so on, till convergence is reached.

At each iteration, the approximation becomes closer to the true function.
Of course, it is very possible that for a generator, its 
assigned power does not change from one iteration to another and so, the approximation will stay the same. This is actually desired because it
means we add less complexity than we might have expected.

\begin{figure}
\centering
\begin{tikzpicture}
\begin{axis} [
	axis x line=bottom,
	axis y line=left,
	xlabel=$x$, ylabel=$y$,
	xtick={-3.14159, -1.5708, 0, 1.5708, 3.14159},
	xticklabels={$-\pi$, $-\frac{\pi}{2}$, 0, $\frac{\pi}{2}$, $\pi$},
	extra x ticks={1},
	extra x tick labels={1}
	]
]
  \addplot[id=4a,domain=-4:4,dotted] gnuplot[samples=250] {abs(sin(x))};
  \addplot[id=4b,domain=-4:4,solid] gnuplot[samples=250]
  						 {0.6366*abs(asin(sin(x)))};
  \addplot[id=4c,domain=-4:4,dashed] gnuplot[samples=250]
  						 {-0.2819*(abs(abs(asin(sin(x)))-1)-1)  						
  						  +0.5596*abs(asin(sin(x)))};
  \draw[dotted] (axis cs:1,0) -- (axis cs:1,0.8415);
\end{axis}
\end{tikzpicture}

\caption{Plot of the two first approximations of the sine term for a given 
         generator. We assume that the solution for this unit in the first
         iteration was 1.  }
\label{approx3}
\end{figure}
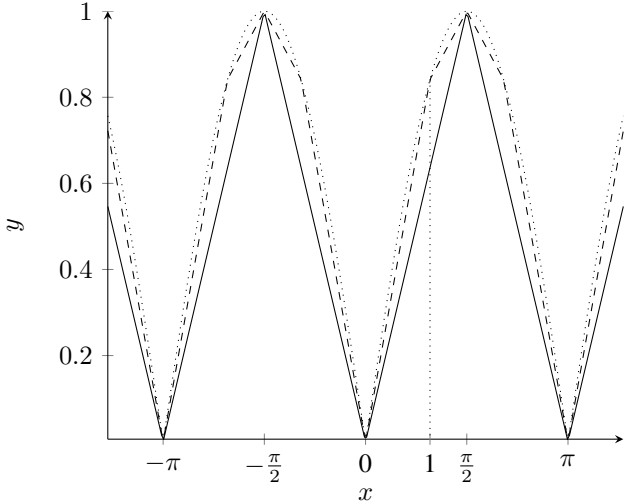

This algorithm can be formalized for our specific needs as follows, where $\mathcal{X}_i$ is the set of break points for cost function
$i$, $m_i$ the number of segments for its approximation, $X_{i,j}$ the $j$th element of $\mathcal{X}_i$ in ascending order, $\alpha_{ij}$
and $\beta_{ij}$ the coefficients of line segment $j$ of cost function $i$, $\delta$ the change in optimal value as a proxy to measure 
convergence, and $\epsilon$ a given positive parameter.

\begin{algorithm}
\caption{Adaptive piecewise-quadratic under-approximation}
\label{al:adaptive}
\begin{algorithmic}[5]
\FOR{$i=1,\dots,n$}
	\STATE $\mathcal{X}_i \gets \{0,\pi/2\}$
	\STATE $m_i \gets 1$
\ENDFOR
\REPEAT 
	\FOR{$i=1,\dots,n$}
		\FOR{$j=1,\dots,m_i$}
			\STATE $\alpha_{i,j} \gets (\sin(X_{i,j+1})-\sin(X_{i,j}))/(X_{i,j+1}-X_{i,j})$
			\STATE $\beta_{i,j} \gets \sin(X_{i,j})-\alpha_{i,j} X_{i,j}$
		\ENDFOR
	\ENDFOR
	\STATE $(\hat{\vec p},\hat{g}) \gets$ solve problem~\eqref{eq:sys_algo}, where $\hat{\vec p}$ denotes the optimal $\vec p$ and $\hat{g}$ the optimal value of the surrogate cost function;
	\STATE $\delta \gets f(\hat{\vec p}) - \hat{g}$ 
	\FOR{$i=1,\dots,n$}
		\STATE $\mathcal{X}_i \gets \mathcal{X}_i \cup \{t_i\}$  \label{step:new-breakpoint}
		\STATE $m_i \gets \#\mathcal{X}_i-1$
	\ENDFOR
\UNTIL{$\delta < \epsilon$}
\end{algorithmic}
\end{algorithm}

We now analyze the convergence of Algorithm~\ref{al:adaptive}.

\begin{theorem}[convergence to the global minimum]
Let $f^*$ denote the optimal value of the cost function $f$~\eqref{eq:f} of the ELDP~\eqref{eq:eldp}. For $m=0,1,\dots$, let $g_m$ denote the piecewise-quadratic surrogate cost function used by Algorithm~\ref{al:adaptive} at iteration $m$, and let $\vec p_m\in\mathbb{R}^n$ denote the power outputs produced by Algorithm~\ref{al:adaptive} at iteration $m$ by solving problem~\eqref{eq:sys_algo}.   
Recall that $g_m(\vec p_m) = \min g_m$ and denote it by $g^*_m$.
Then $\lim_{m\to\infty}g^*_m = \lim_{m\to\infty}f(\vec p_m) = f^*$, and every limit point of $(\vec p_m)_{m\in\mathbb{N}}$ is a global solution of the ELDP.  
\end{theorem}
\begin{IEEEproof}
We first show that $f$ and $g_m$, $m=0,1,\dots$, are Lipschitz continuous on the ELDP feasible set with a common Lipschitz constant $K$. Indeed, for every $i$, the cost function $f_i$~\eqref{eq:fi} satisfies the Lipschitz property $|f_i(p_i+\Delta)-f_i(p_i)| \leq (2a_i p^{\mathrm{max}}_i + b_i + d_ie_i)\Delta$ for all $p_i$ and $p_i+\Delta$ that satisfy the generator capacity constraints~\eqref{eq:box}. The ELDP cost function $f$~\eqref{eq:f}, being the sum of Lipschitz continuous functions, is thus Lipschitz continuous on the ELDP feasible set, with a constant $K=\sum_{i=1}^n 2a_i p^{\mathrm{max}}_i + b_i + d_ie_i$. 
Since $g_m$ is obtained by replacing the rectified sines of $f$ by chords, it follows that $2a_i p^{\mathrm{max}}_i + b_i + d_ie_i$ is still a Lipschitz constant for the contribution of generator $i$ to $g_m$, and hence that $K$ is also a Lipschitz constant for $g_m$.

Since building $g_{m+1}$ from $g_m$ consists of inserting, for each generator, one new breakpoint for the piecewise-linear under-approximation of the piecewise-concave rectified sine, it follows that, for every $m$, 
\begin{equation}
\label{eq:seqbnd}
g_m \leq g_{m+1} \leq f.
\end{equation}
Therefore $(g^*_m)_{m\in \mathbb{N}}$ is a nondecreasing sequence bounded by $f^*$.
Thus $(g^*_m)_{m\in \mathbb{N}}$ converges and $\lim_{m\to\infty}g^*_m\leq f^*$.

We show that $\lim_{m\to\infty}g^*_m = f^*$. By contradiction, assume that $\lim_{m\rightarrow \infty} g^*_m = f^*-\varepsilon$ with $\varepsilon>0$. 
Since $(\vec p_m)_{m\in\mathbb{N}}$ is bounded in view of the generator capacity constraints~\eqref{eq:box}, there
exists a subsequence $(\vec p_{m_k})_{k\in\mathbb{N}}$ that converges.
We then have the following inequalities which we justify hereafter:
\begin{align}
\|\vec p_{m_k}-\vec p_{m_{k+1}}\| &\geq \frac{1}{K} [g_{m_{k+1}}(\vec p_{m_k}) - g_{m_{k+1}}(\vec p_{m_{k+1}}) ]  \label{eq:bnd1}
\\ &\geq \frac{1}{K} [f^* - g_{m_{k+1}}(\vec p_{m_{k+1}}) ]  \label{eq:bnd2}
\\ &\geq \frac{1}{K} [f^* - (f^*-\varepsilon) ]  \label{eq:bnd3} \\ 
&\geq \frac{\varepsilon}{K}, \nonumber
\end{align}
a contradiction since $(\vec p_{m_k})_{k\in\mathbb{N}}$ converges. Inequality~\eqref{eq:bnd1} states Lipschitz continuity of $g_{m_{k+1}}$ with constant $K$. Inequality~\eqref{eq:bnd2} follows from $g_{m_{k+1}}(\vec p_{m_{k}}) = f(\vec p_{m_k}) \geq f^*$; indeed, by construction, $\vec p_{m_k}$ is a breakpoint of $g_{m_k+1}$ and all subsequent surrogate cost functions. Finally,~\eqref{eq:bnd3} follows from $g_{m_{k+1}}(\vec p_{m_{k+1}}) = g^*_{m_{k+1}} \leq f^*-\varepsilon$. 

We now show that $\lim_{m\to\infty}f(\vec p_m)=f^*$. By contradiction, suppose not. Then there is an infinite subsequence $(\vec p_{m_k})_{k\in\mathbb{N}}$ and $\varepsilon>0$ such that $f(\vec p_{m_k}) \geq f^*+\epsilon$. We assume w.l.o.g.\ that $(\vec p_{m_k})_{k\in\mathbb{N}}$ converges; if not, we extract a sub-subsequence that does. By the triangle inequality, we have
\begin{align*}
|f(\vec p_{m_k}) - f^*| \leq & |f(\vec p_{m_k}) - g_{m_{k+1}}(\vec p_{m_k})|
\\ + & |g_{m_{k+1}}(\vec p_{m_k}) - g_{m_{k+1}}(\vec p_{m_{k+1}})|
\\ + & |g_{m_{k+1}}(\vec p_{m_{k+1}}) - f^*|.
\end{align*}
The first term of the bound is zero by construction of the surrogate functions. The second term goes to zero as $k\to\infty$ in view of the common Lipschitz constant $K$ and the convergence of $(\vec p_{m_k})_{k\in\mathbb{N}}$. The third term goes to zero as $k\to\infty$ since $g_{m_{k+1}}(\vec p_{m_{k+1}})= g^*_{m_{k+1}}$ and $\lim_{m\to\infty}g^*_m = f^*$. Hence $|f(\vec p_{m_k}) - f^*|$ goes to zero, a contradiction.

Finally, let $(\vec p_{m_k})_{k\in\mathbb{N}}$ be a convergent subsequence and let $\vec p_{m_\infty}$ denote its limit. By continuity of $f$, we have that $f(\vec p_{m_\infty}) = \lim_{k\to\infty} f(\vec p_{m_k}) = f^*$. 
\end{IEEEproof}
Note that, in view of the breakpoint insertion procedure (step~\ref{step:new-breakpoint} of Algorithm~\ref{al:adaptive}), it generally does not hold that $g_m$ converges to $f$ pointwise; otherwise the proof above could have been more direct. 

\section{Numerical results}
\label{sec:results}

After having built such a model, we can hand it over to a solver that can
handle MIQP. Here we will use the Gurobi solver~\cite{Gurobi}. To study the efficiency
of our method, we will test it on the most common test cases in the
literature : a 3-units setting with a demand of 850 MW~\cite{Walters93} (I), a 13-units setting
with a demand of 1800 MW~\cite{Sinha03} (IIa) and 2520 MW~\cite{Wong94} (IIb), and a 40-units setting with a demand of 10,500 MW~\cite{Sinha03} (III). 
Using model \eqref{eq:simple} and these datasets, and feeding them to
Gurobi, we get the solutions given in tables~\ref{tab:simple-3}, \ref{tab:simple-13}, \ref{tab:simple-40}. The hardware used is a PC
with a Intel Core 2 Duo P8600 CPU (two cores at 2,4 Ghz) and 3 GB of RAM,
on GNU/Linux.

\begin{table}[!t]
\renewcommand{\arraystretch}{1.3}

\caption{Found solution for the 3-units study case (I) using the simple model \eqref{eq:simple}}
\label{tab:simple-3}

\centering
\begin{tabular}{>{$}c<{$} c}
\hline
\text{\bfseries Unit} & \bfseries Power (MW) \\ \hline
p_1 & 300.267 \\
p_2 & 400.000 \\ 
p_3 & 149.733 \\ \hline
\text{\bfseries Total cost (\$/h)} & 8234.07 \\
\text{\bfseries Best in lit. (\$/h)} & 8234.07 \\
\text{\bfseries Real time (s)} & 0.013 \\
\text{\bfseries CPU time (s)} & 0.012 \\ \hline
\end{tabular}

\end{table}

\begin{table}[!t]
\renewcommand{\arraystretch}{1.3}

\caption{Found solution for the 13-units study case (IIa) using the simple model \eqref{eq:simple}}
\label{tab:simple-13}

\centering
\begin{tabular}{cr}
\hline
\bfseries Unit & \bfseries Power (MW) \\ \hline
$p_1$  & 628.319 \\ 
$p_2$  & 222.749 \\
$p_3$  & 149.600 \\ 
$p_4$  & 109.867 \\ 
$p_5$  &  60.000 \\ 
$p_6$  & 109.867 \\ 
$p_7$  & 109.867 \\ 
$p_8$  & 109.867 \\ 
$p_9$  & 109.867 \\ 
$p_{10}$  & 40.000 \\ 
$p_{11}$  & 40.000 \\ 
$p_{12}$  & 55.000 \\ 
$p_{13}$  & 55.000 \\ \hline
\bfseries Total cost (\$/h) & 17963.83 \\ 
\bfseries Best in lit. (\$/h) &  17963.83 \\ 
\bfseries Real time (s) & 0.381 \\ 
\bfseries CPU time (s) & 0.708 \\
\hline
\end{tabular}
\end{table}

\begin{table}
\caption{Found solution for the 13-unit case (IIb) using the simple model \eqref{eq:simple}}
\label{tab:simple2b}
\centering
\begin{tabular}{cc}
\hline
\bfseries Unit & \bfseries Power (MW) \\
\hline
$p_1$  & 628.319 \\ 
$p_2$  & 299.199 \\
$p_3$  & 299.199 \\ 
$p_4$  & 159.733 \\ 
$p_5$  & 159.733 \\ 
$p_6$  & 159.733 \\ 
$p_7$  & 159.733 \\ 
$p_8$  & 159.733 \\ 
$p_9$  & 159.733 \\ 
$p_{10}$  & 77.400 \\ 
$p_{11}$  & 77.400 \\ 
$p_{12}$  & 90.042 \\ 
$p_{13}$  & 90.042 \\ 
\hline
\bfseries Total cost (\$/h) & 24170.66 \\ 
\bfseries Best in lit. (\$/h) &  24169.92 \\ 
\bfseries Real time (s) & 0.054 \\ 
\bfseries CPU time (s) & 0.084 \\
\hline
\end{tabular}
\end{table}

\begin{table}[!t]
\renewcommand{\arraystretch}{1.3}

\caption{Found solution for the 40-units study case (III) using the simple model \eqref{eq:simple}}
\label{tab:simple-40}

\centering
    \begin{tabular}{cr|cr}
    \hline
\bfseries Unit & \bfseries Power (MW) & \bfseries Unit & \bfseries Power (MW)\\ \hline
    $p_1$ & 110.800  &          $p_{21}$ & 523.279 \\
    $p_2$ & 110.800  &          $p_{22}$ & 523.279 \\
    $p_3$ & 97.400   &          $p_{23}$ & 523.279 \\
    $p_4$ & 179.733 &           $p_{24}$ & 523.279 \\
    $p_5$ & 90.279  &           $p_{25}$ & 523.279 \\
    $p_6$ & 140.000  &            $p_{26}$ & 523.279 \\
    $p_7$ & 259.600  &          $p_{27}$ & 10.000  \\
    $p_8$ & 284.600  &          $p_{28}$ & 10.000  \\
    $p_9$ & 284.600  &          $p_{29}$ & 10.000  \\
    $p_{10}$ & 130.000  &         $p_{30}$ & 90.279  \\
    $p_{11}$ & 168.800  &       $p_{31}$ & 190.000 \\
    $p_{12}$ & 168.800  &       $p_{32}$ & 190.000 \\
    $p_{13}$ & 214.760  &       $p_{33}$ & 190.000 \\
    $p_{14}$ & 394.279 &        $p_{34}$ & 164.800 \\
    $p_{15}$ & 394.279 &        $p_{35}$ & 164.800 \\
    $p_{16}$ & 304.520 &        $p_{36}$ & 164.800 \\
    $p_{17}$ & 489.279 &        $p_{37}$ & 110.000 \\
    $p_{18}$ & 489.279 &        $p_{38}$ & 110.000 \\
    $p_{19}$ & 511.279  &       $p_{39}$ & 110.000 \\
    $p_{20}$ & 511.279  &       $p_{40}$ & 511.279 \\
    \hline
\multicolumn{2}{c}{\bfseries Total cost (\$/h)} & \multicolumn{2}{c}{121415.31} \\
\multicolumn{2}{c}{\bfseries Best in lit. (\$/h)} & \multicolumn{2}{c}{121412.54} \\ 
\multicolumn{2}{c}{\bfseries Real time (s)} & \multicolumn{2}{c}{0.116} \\ 
\multicolumn{2}{c}{\bfseries CPU time (s)} & \multicolumn{2}{c}{0.192} \\
\hline
\end{tabular}    
\end{table}

If we compare these results to those found in the literature~\cite{Srinivasa13}, 
we can see
that they are close to the best solutions found to date. 
In fact, for I and IIa, we do get the best results, while for IIb and III,
it is less than 0.005\% worse. Furthermore, it is achieved on a deterministic basis, without the uncertainty and irreproducibility of
the commonly used heuristics.

For the model of Section~\ref{sec:finer_solver}, table~\ref{solution} shows the solution found for a choice of parameter $\theta_1 = 0.35\pi$ and $\theta_2=0.47\pi$. 

The resulting optimal value is slightly better than what we found earlier. The real time needed is 0.493 s for a CPU time of 0.902 s.

\begin{table}[!t]
\renewcommand{\arraystretch}{1.3}

\caption{Solution for 40-units case study (III) using the model of Section~\ref{sec:finer_solver}.}
\label{solution}

\centering

\begin{tabular}{cr|cr}
\hline
\bfseries Unit & \bfseries Power (MW) & \bfseries Unit & \bfseries Power (MW)\\ \hline
 $p_1$    & 110.800 & $p_{21}$ & 523.279 \\                    
 $p_2$    & 110.800 & $p_{22}$ & 523.279 \\                     
 $p_3$    &  97.400 & $p_{23}$ & 523.279 \\                     
 $p_4$    & 179.733 & $p_{24}$ & 523.279 \\                     
 $p_5$    &  87.800 & $p_{25}$ & 523.279 \\    
 $p_6$    & 140.000 & $p_{26}$ & 523.279 \\    
 $p_7$    & 259.600 & $p_{27}$ &  10.000 \\    
 $p_8$    & 284.600 & $p_{28}$ &  10.000 \\    
 $p_9 $   & 284.600 & $p_{29}$ &  10.000 \\
 $p_{10}$ & 130.000 & $p_{30}$ &  87.800 \\
 $p_{11}$ &  94.000 & $p_{31}$ & 190.000 \\
 $p_{12}$ &  94.000 & $p_{32}$ & 190.000 \\
 $p_{13}$ & 214.760 & $p_{33}$ & 190.000 \\
 $p_{14}$ & 394.279 & $p_{34}$ & 164.800 \\
 $p_{15}$ & 394.279 & $p_{35}$ & 200.000 \\
 $p_{16}$ & 394.279 & $p_{36}$ & 194.398 \\
 $p_{17}$ & 489.279 & $p_{37}$ & 110.000 \\
 $p_{18}$ & 489.279 & $p_{38}$ & 110.000 \\
 $p_{19}$ & 511.279 & $p_{39}$ & 110.000 \\
 $p_{20}$ & 511.279 & $p_{40}$ & 511.279 \\ \hline
\multicolumn{2}{c}{\bfseries Total cost (\$/h)} & \multicolumn{2}{c}{121412.54} \\
\multicolumn{2}{c}{\bfseries Best in lit. (\$/h)} & \multicolumn{2}{c}{121412.54} \\ 
\multicolumn{2}{c}{\bfseries Real time (s)} & \multicolumn{2}{c}{0.493} \\ 
\multicolumn{2}{c}{\bfseries CPU time (s)} & \multicolumn{2}{c}{0.902} \\
\hline

\end{tabular}
\end{table}

We also tested the method of Section~\ref{sec:global_convergence} on the different study cases. It so happens
that it finds exactly the same solutions (except for a few swaps in equivalent
generators) as the method of Section~\ref{sec:finer_solver} for case I, IIa and III
thereby proving their optimality. The exception is case IIb, for which
a slightly better solution is found (table \ref{tab:global2b}).
Note that the algorithm of Section~\ref{sec:global_convergence} always needs more time 
since it takes at least two iterations to stop: one to get a solution and
another to prove the optimality.

\begin{table}
\caption{Global solution for case IIb found with the algorithm of Section~\ref{sec:global_convergence}}
\label{tab:global2b}

\centering
\begin{tabular}{cc}
\hline
\bfseries Unit & \bfseries Power (MW) \\
\hline
$p_1$  & 628.319 \\ 
$p_2$  & 299.199 \\
$p_3$  & 299.199 \\ 
$p_4$  & 159.733 \\ 
$p_5$  & 159.733 \\ 
$p_6$  & 159.733 \\ 
$p_7$  & 159.733 \\ 
$p_8$  & 159.733 \\ 
$p_9$  & 159.733 \\ 
$p_{10}$  & 77.400 \\ 
$p_{11}$  & 77.400 \\ 
$p_{12}$  & 87.684 \\ 
$p_{13}$  & 92.400 \\ 
\hline
\bfseries Total cost (\$/h) & 24169.92 \\ 
\bfseries Best in lit. (\$/h) &  24169.92 \\ 
\bfseries Real time (s) & 0.589 \\ 
\bfseries CPU time (s) & 0.592 \\
\hline
\end{tabular}
\end{table}

\section{Conclusion}
\label{sec:conclusion}

This article concerns a piecewise quadratic under-approximation of the ELDP cost function that takes into account the valve-point effect, thereby providing a means to solve it globally with the MIQP method, despite the nonsmoothness and the nonconvexity of the original cost function. Furthermore, a convergence analysis is presented to show that, under mild assumptions, this strategy guarantees the global minimizer in the static ELDP context. In order to support our claim, the minimization results are furnished for the datasets corresponding to a 3-, two instances of 13-, and a 40-generator setting, wherein the transmission losses are omitted. 
Interestingly enough, one the one hand, in accordance with the global convergence guarantee, the outcome of the cost minimization by the method of Section~\ref{sec:global_convergence} is never surpassed by the hybrid methods, and on the other hand, the computational times are quite impressive.

While these widely used datasets in the ELDP literature enable us to investigate how our approach compares with the state-of-the-art methods, the same framework is applicable for scenarios that do not ignore losses. 
The only modification needed is to add the loss term in the constraint
to take it into account and to relax the equality into an inequality.
The reason behind the relaxation is to make the constraint convex and
the model easy enough to solve. This does not change the solution of the
problem as long as the objective function is monotonically increasing.

Moreover, these methods can be used for other cost functions, where
the valve-point effect would be modeled differently. As long as the
valve-point effect is represented by a periodic piecewise-concave function, only
minor adaptation would be needed.

\bibliographystyle{IEEEtran}
\bibliography{IEEEabrv,ELDP_bibliography}

\end{document}